\documentstyle[11pt]{article}
\def\fa{f_{\alpha}}

\def\fia{\phi_{\alpha}}
\def\pia{\psi_{\alpha}}

\def\ga{\Gamma_{\alpha}}
\def\da{\Delta_{\alpha}}
\def\gb{\G_{\beta}}

\def\gab{g_{\alpha\beta}}

\def\va{V_{\alpha}}
\def\wa{W _{\alpha}}
\def\ut{\tilde{U}}
\def\wt{\tilde{W}}

\def\vta{\vt_{\alpha}}
\def\wta{\wt_{\alpha}}

\def\vtb{\vt_{\beta}}
\def\oa{\omega_\alpha}

\def\vb{V_{\beta}}
\def\fib{\phi_{\beta}}

\def\G{\Gamma}

\def\ut{\tilde{U}}

\def\vt{\tilde{V}}

\def\ft{\tilde{f}}

\def\ot{\tilde{\omega}}

\def\D{\Delta}

\def\c0{{\cal C}_0}

\def\C{\mbox{\bbb{C}}}
\def\N{\mbox{\bbb{N}}}

\def\Q{\mbox{\bbb{Q}}}
\def\R{\mbox{\bbb{R}}}
\def\Z{\mbox{\bbb{Z}}}

\def\ctwo{\C^2}
\def\cthree{\C^3}

\def\cd{\C^d}
\def\ck{\C^k}

\def\rtwo{\R^2}

\def\rd{\R^d}

\def\rdu{\R^*}
\def\rtwodu{(\rtwo)^*}

\def\rk{\R^k}

\def\Zd{\Z^d}

\def\vz{\underline{z}}

\def\z1s{|z_1|^2}
\def\ztwos{|z_2|^2}

\def\et1{e^{2\pi i\theta_1}}

\def\ed{e_1,\ldots,e_d}
\def\ld{\lambda_1,\ldots,\lambda_d}

\def\xd{X_1,\ldots,X_d}

\def\b{\mbox{\frak b}}
\def\d{\mbox{\frak d}}

\def\n{\mbox{\frak n}}
\def\t{\mbox{\frak t}}
\def\ddu{\d^*}

\def\X{\mbox{X}}

\def\ck{\C^k}
\def\zd{(z_1,\cdots,z_d)}
\def\lorw{\longrightarrow}

\def\vsh{V^{\#}}
\def\SC{\mbox{\sbbb{C}}}

\def\Dc{D_{\SC}}
\def\dc{\d_{\SC}}
\def\bc{\b_{\SC}}
\def\Tc{T_{\SC}}
\def\Tdc{T^d_{\SC}}

\def\Nc{N_{\SC}}
\def\nc{\n_{\SC}}
\def\tc{\t_{\SC}}
\def\cdf{\C^d_F}
\def\cdd{\C^d_{\Delta}}
\def\cddtwo{\C^2_{\Delta}}
\def\cddth{\C^3_{\Delta}}
\def\cddfv{\C^5_{\Delta}}

\def\vtn{\tilde{V}_{\sc n}}
\def\vts{\tilde{V}_{\sc s}}

\def\vn{V_{\sc n}}
\def\vs{V_{\sc s}}
\def\pic{\pi_{\SC}}
\def\Pic{\Pi_{\SC}}

\def\squo{\Psi^{-1}(0)/N}
\def\cquo{\cdd/\Nc}
\def\zset{\Psi^{-1}(0)}

\def\smuc{\mbox{\frak{s}}^{\mu}_{\SC}}

\def\Smuc{S^{\mu}_{\SC}}

\def\vmu{V_{\mu}}
\def\vcmu{\widehat{V}_{\mu}}
\def\vtmu{\tilde{V}_{\mu}}
\def\cdmu{\cd_{\mu}}
\def\vnu{V_{\nu}}

\def\vtnu{\tilde{V}_{\nu}}

\def\wmu{W_{\mu}}
\def\wtmu{\tilde{W}_{\mu}}
\def\wfmu{W^{\#}_{\mu}}
\def\wnu{W_{\nu}}
\def\wtnu{\tilde{W}_{\nu}}
\def\wfnu{W^{\#}_{\nu}}
\def\fimu{\phi_{\mu}}
\def\finu{\phi_{\nu}}
\def\gimunu{g_{\mu\nu}}
\def\gimunus{\gimunu^{\#}}

\def\vw{\underline{w}}
\def\vx{\underline{x}}
\def\vy{\underline{y}}
\def\veta{\underline{\eta}}
\def\vzeta{\underline{\zeta}}
\def\vtone{\tilde{V}_1}
\def\vttwo{\tilde{V}_2}
\def\vtthree{\tilde{V}_3}

\def\s{\mbox{\frak{s}}}
\def\un{U_{\sc n}}
\def\us{U_{\sc s}}

\def\utn{\tilde{U}_{\sc n}}
\def\uts{\tilde{U}_{\sc s}}
\def\chin{\tilde{\chi}_{\sc n}}
\def\chis{\tilde{\chi}_{\sc s}}
\def\wtn{W^{\#}_{\sc n}}
\def\wts{W^{\#}_{\sc s}}

\def\umu{U_{\mu}}
\def\ucmu{\widehat{U}_{\mu}}
\def\utmu{\tilde{U}_{\mu}}

\def\chitmu{\tilde{\chi}_{\mu}}

\def\vu{\underline u}
\def\gmu{\Gamma_{\mu}}
\def\s{\mbox{\frak{s}}}

\newtheorem{thm}{Theorem}[section]
\newtheorem{prop}[thm]{Proposition}
\newtheorem{lemma}[thm]{Lemma}

\newtheorem{defn}[thm]{Definition}
\newtheorem{remark}[thm]{Remark}

\newtheorem{ex}[thm]{Example}

\newcommand{\qed}{\nolinebreak\hfill{$\Box$}\par\vspace{0.5\parskip}}

\newcommand{\proof}{\mbox{\bf Proof.\ \ }}
\newfont{\frak}{eufm10 scaled\magstep1}
\newfont{\sfrak}{eufm8 scaled\magstep1}
\newfont{\bbb}{msbm10 scaled\magstephalf}
\newfont{\sbbb}{msbm7 scaled\magstephalf}

\newcounter{sect}

\title{\sc Generalized Toric Varieties for Simple Non-Rational Convex Polytopes}
\author{\sc Fiammetta Battaglia and Elisa Prato}
\date{}
\begin{document}
\maketitle
\begin{abstract}
We call {\em complex quasifold} of dimension $k$ a space that is locally isomorphic
to the quotient of an open subset of the space $\ck$ by the holomorphic action of a
discrete group; the analogue of a complex torus in this setting is called a {\em
complex quasitorus}. We associate to each simple polytope, rational or not, a family
of complex quasifolds having same dimension as the polytope, each containing a dense
open orbit for the action of a suitable complex quasitorus. We show that each of
these spaces $M$ is diffeomorphic to one of the symplectic quasifolds defined in
\cite{p}, and that the induced symplectic structure is compatible with the complex
one, thus defining on $M$ the structure of a K\"ahler quasifold. These spaces may be
viewed as a generalization of the toric varieties that are usually associated to
those simple convex polytopes that are rational.
\end{abstract}
\section*{Introduction}
Consider a vector space $\d$ of dimension $n$. To each simple convex polytope
$\D\subset\ddu$ that is rational with respect to a lattice in $\d$ there corresponds
a toric variety with at worst quotient singularities. What happens in the case that
the simple convex polytope is no longer rational? To answer this question we consider
a special class of spaces, called {\em quasifolds}, which were first introduced by
one of the authors in \cite{p}. A quasifold is not necessarily a Hausdorff space: it
is locally modeled by orbit spaces of the action of discrete, possibly infinite,
groups on open subsets of $\rk$. A {\em quasitorus}, on the other hand, is the
natural replacement of a torus in this geometry.

In this article we define the notions of complex quasifold and complex quasitorus and
we associate to each simple convex polytope $\Delta\subset\ddu$ a family of compact
complex quasifolds of dimension $n$, each endowed with the holomorphic action of a
complex quasitorus $\Dc$ having a dense open orbit. Our construction is explicit:
each space $M$ is the topological quotient of a suitable open subset of $\cd$ by the
action of a suitable subgroup $\Nc\subset\Tdc=\cd/\Zd$, and $\Dc$ is isomorphic to
$\Tdc/\Nc$, $d$ being the number of facets of the polytope. We show that $M$ is a
complex quasifold by covering it with mutually compatible local models of the type
$\C^n$ modulo the action of a discrete group. If the polytope is rational our
procedure matches the standard one for constructing toric varieties as quotients that
is described in Appendix~1 of Guillemin's book \cite{vg}.

It is shown in \cite{p} that to the same simple convex polytope $\D$ one can also
associate a family of symplectic quasifolds of dimension $2n$, each endowed with the
effective Hamiltonian action of a quasitorus $D$, having the property that the image
of the corresponding moment mapping is exactly $\D$. The construction extends the one
given by Delzant \cite{d} in the smooth case: it is explicit and uses the symplectic
quotient operation. In the last section of this paper we show that each complex
quotient $M$ is equivariantly diffeomorphic to one of these symplectic quotients,
that the induced symplectic structure is compatible with the complex one, and thus
defines on $M$ the structure of a K\"ahler quasifold; of course $\Dc$ here is the
complexification of the corresponding quasitorus $D$.

For these reasons these spaces may well be thought of as a natural generalization of
toric varieties for simple convex polytopes that are not rational.

\section{Complex quasifolds and complex quasitori}\label{folds}
This section is devoted to defining complex quasifolds and their geometry. We will
not repeat the remarks and results that are in common with the real case, for which
we refer the reader to the article \cite{p}. The local model for complex quasifolds
is a complex manifold acted on holomorphically by a discrete group.
\begin{defn}[Complex model]\label{model}{\rm
    Let $\vt$ be a connected, simply connected complex manifold of
    dimension $k$ and let $\G$ be a discrete group acting on $\vt$
    holomorphically so that the set of points, $\vt_0$, where the action is free,
    is connected and dense. Consider the space of orbits, $\vt/\G$, of the action of the group
    $\G$ on the manifold $\vt$, endowed with the quotient topology, and the
    canonical projection $p\;\colon\;\vt\rightarrow \vt/\G$.
    A {\em complex model} of dimension $k$ is the triple $(\vt/\G,p,\vt)$, shortly
    $\vt/\G$.}
\end{defn}
\begin{remark}\label{simplyc}{\rm
We remark that the assumption in Definition~\ref{model} that the manifold $\vt$ be
simply connected is not as strong as one may think. Consider the triple
$(\vt/\G,p,\vt)$ as defined above but assume that the manifold $\vt$ is not simply
connected; consider its universal cover, $\pi\,\colon\,\vsh\rightarrow\vt$, and its
fundamental group, $\Pi$. The manifold $\vsh$ is connected and simply connected, the
mapping $\pi$ is holomorphic, the discrete group $\Pi$ acts holomorphically, freely
and properly on the manifold $\vsh$ and $\vt=\vsh/\Pi$. Consider the extension of the
group $\G$ by the group $\Pi$, $1\longrightarrow \Pi\longrightarrow\Lambda
\longrightarrow\G\longrightarrow 1$, defined as follows
$$\Lambda=\left\{\;\lambda\in\mbox{Diff}(\vsh)\;|\;\exists\;
\gamma\in\Gamma\;\mbox{s. t.}\;\pi(\lambda(u^{\#}))=\gamma\cdot
\pi(u^{\#})\;\forall\; u^{\#}\in\vsh\;\right\}.$$ It is easy to verify that $\Lambda$
is a discrete group, that it acts on the manifold $\vsh$ according to the assumptions
of Definition~\ref{model} and that $\vt/\G=\vsh/\Lambda$.}\end{remark}
\begin{defn}[Submodel]{\rm
Consider a model $(\vt/\G,p,\vt)$ and let $W$ be an open subset of $\vt/\G$. We will
say that $W$ is a submodel of $(\vt/\G,p,\vt)$, if $(p^{-1}(W),p,W)$ defines a model
according to Remark~\ref{simplyc}.}
\end{defn}
\begin{defn}[Holomorphic mapping, biholomorphism of models]{\rm
Given two models $(\vt/\G, p, \vt)$ and $(\wt/\D,q,\wt)$, a mapping
$f\,\colon\,\vt/\G\longrightarrow\wt/\D$ is said to be {\em holomorphic} if there
exists a holomorphic mapping $\ft\,\colon\,\vt\longrightarrow\wt$ such that $q\circ
\ft= f\circ p$; we will say that $\ft$ is a {\em lift} of $f$. We will say that the
holomorphic mapping $f$ is a {\em biholomorphism of models} if it is bijective and if
the lift $\ft$ is a biholomorphism.}
\end{defn}
If the mapping $\ft$ is a lift of a holomorphic mapping of models
$f\,\colon\,\ut/\G\longrightarrow \vt/\D$ so are the mappings
$\ft^{\gamma}(-)=\ft(\gamma\cdot -)$, for all elements $\gamma$ in $\G$ and
$^{\delta}\ft(-)=\delta\cdot\ft(-)$, for all elements $\delta$ in $\Delta$. If the
mapping $f$ is a biholomorphism, then these are the only other possible lifts and the
groups $\G$ and $\D$ are isomorphic; the proof goes exactly as in the real case, see
\cite[orange and green lemmas]{p}.

Geometric objects on a model $\vt/\G$ are defined by $\G$-invariant geometric objects
on the manifold $\vt$. For example:
\begin{defn}[Differential form on a model]\label{locform}{\rm
A {\em differential form of degree $h$}, $\omega$, on a model $\vt/\G$ is the
assignment of a $\G$-invariant differential form of degree $h$, $\ot$, on the complex
manifold $\vt$.}
\end{defn}
\begin{defn}[K\"ahler form on a model]{\rm A {\em K\"ahler form}
on a model $\vt/\G$ is a differential form, $\omega$, such that $\ot$ (see
Definition~\ref{locform}) is K\"ahler.}
\end{defn}
Complex quasifolds are obtained by gluing together the models in the appropriate way:
\begin{defn}[Complex quasifold]
  {\rm A dimension $k$ {\em complex quasifold structure} on a topological space
    $M$ is the assignment of an {\em atlas}, or collection of {\em charts},
   ${\cal A}= \{\,(\va,\fia,\vta/\ga)\,|\,\alpha\in A\,\}$ having the following properties:
\begin{enumerate}
\item The collection $\{\,\va\,|\,\alpha\in A\,\}$ is a cover of $M$.
\item For each index $\alpha$ in $A$, the set $\va$ is open, the space
$\vta/\ga$ defines a model, where the set $\vta$ is an open, connected, and simply
connected subset of the space $\ck$, and the mapping $\fia$ is a homeomorphism of the
space $\vta/\ga$ onto the set $\va$.
\item For all indices $\alpha, \beta$ in $A$ such that
  $\va\cap\vb\neq\emptyset$, the sets $\fia^{-1}(\va\cap\vb)$ and
  $\fib^{-1}(\va\cap\vb)$ are submodels of $\vta/\ga$ and $\vtb/\gb$
  respectively and the mapping
  $$\gab=\fib^{-1}\circ\fia\,\colon\fia^{-1}(\va\cap\vb)
  \longrightarrow\fib^{-1}(\va\cap\vb)$$
  is a biholomorphism of models. We will then say that the mapping
  $\gab$ is a {\em change of charts} and that the corresponding charts are
  {\em compatible}.
\item The atlas $\cal A$ is maximal, that is: if the triple
$(V,\phi,\vt/\G)$ satisfies property 2. and is compatible with all the charts in
$\cal A$, then $(V,\phi,\vt/\G)$ belongs to $\cal A$.
\end{enumerate}
We will say that a space $M$ with a complex quasifold structure is a {\em complex
quasifold}. }
\end{defn}
A complex quasifold of dimension $k$ has an underlying structure of real quasifold of
dimension $2k$.
\begin{remark}{\rm To each point $m\in M$ there corresponds a discrete group $\G_m$
defined as follows: take a chart $(\va,\fia,\vta/\ga)$ around $m$, then $\G_m$ is the
isotropy group of $\ga$ at any point $\tilde{v}\in \vt$ which projects down to $m$.
One can check that this definition does not depend on the choice of the
chart.}\end{remark}
\begin{defn}[Holomorphic mapping, biholomorphism]
{\rm Let $M$ and $N$ be two complex quasifolds. A continuous mapping
$f\,\colon\,M\longrightarrow N$ is said to be a {\em holomorphic mapping} if there
exists a chart $(\va,\fia,\vta/\ga)$ around each point $m$ in the space $M$, a chart
$(\wa,\pia,\wta/\da)$ around the point $f(m)$, and a holomorphic mapping of models
$\fa\,\colon\,\vta/\ga\rightarrow\wta/\da$ such that $\pia\circ\fa=f\circ\fia$. If
$f$ is bijective, and if each $\fa$ is a biholomorphism of models, we will say that
$f$ is a {\em biholomorphism}.}
\end{defn}
We remark that two biholomorphic quasifolds are also diffeomorphic with respect to
their underlying real quasifold structure.

Geometric objects on quasifolds are defined as geometric objects on the charts that
respect the changes of charts. For example:
\begin{defn}[Differential form]\label{form}
{\rm A {\em differential form of degree $h$}, $\omega$, on a complex quasifold $M$ is
the assignment of a chart $(\va,\fia,\vta/\ga)$ around each point $m\in M$ and of a
differential form of degree $h$, $\oa$, on the model $\vta/\ga$. We require that
whenever we have two such charts, $(\va,\fia,\vta/\ga)$ and $(\vb,\fib,\vtb/\gb)$,
with the property that $\va\cap\vb\neq\emptyset$, then
$(\gab)^*\omega_{\beta}=\omega_{\alpha}$ for the corresponding change of charts
$\gab$.}
\end{defn}
\begin{defn}[K\"ahler form]\label{kahform}{\rm A {\em K\"ahler
form}, $\omega$, on a complex quasifold $M$ is a differential form, $\omega$, such
that each $\oa$ (see Definition~\ref{form}) is a K\"ahler form on the model
$\vta/\ga$. A {\em K\"ahler structure} on a quasifold $M$ is the assignment of a
K\"ahler form, $\omega$, and we will say that $(M,\omega)$, or simply $M$, is a {\em
K\"ahler quasifold}.}
\end{defn}

The analogue of a torus in the real setting is a quasitorus. Let $\d$ be a vector
space of dimension $n$. We recall from \cite{p} that a quasitorus of dimension $n$
and quasi-Lie algebra $\d$ is the quotient of the space $\d$ by a quasilattice $Q$,
which in turn is the $\Z$-span of a set of spanning vectors in $\d$. Consider now the
complexification of $\d$, the complex vector space $\dc=\d+i\d$; the quasilattice $Q$
acts naturally on $\dc$:
\begin{equation}
\begin{array}{lcccc}
Q&\times&\dc & \lorw & \dc\\ (A&,&X+iY) & \longmapsto & (X+A)+iY.
\end{array}
\label{reallattice}
\end{equation}
\begin{defn}[Complex quasitorus, quasi-Lie algebra, exponential]
{\rm Let $\d$ be a vector space of dimension $n$ and let $Q$ be a quasilattice in
$\d$. Then we call the quotient $\Dc=\dc/Q$ {\em complex quasitorus} of dimension $n$
and {\em quasi-Lie algebra} $\dc$. We call the corresponding projection
$\dc\rightarrow\Dc$ {\em exponential mapping} and we denote it by $\exp$.}
\end{defn}
If the quasilattice $Q$ is a lattice $L$ we obtain the honest complex torus $\dc/L$,
which is the complexification of the torus $\d/L$. The complex quasitorus $\Dc$ is a
quasifold of one chart and may be naturally thought of as the complexification of the
quasitorus $D=\d/Q$. The main result that we will be needing is the following:
\begin{prop}\label{lattice}
Let $T$ be a torus and $N$ a Lie subgroup. Then $\Tc/\Nc$ is a complex quasitorus of
complex dimension $n=\dim T-\dim N$.
\end{prop}
\proof  This proposition has a real analogue, the proof of which is quite similar
(compare with \cite[Proposition~2.5]{p}). Denote by $\n$ the Lie algebra of $N$, and
by $\t$ the Lie algebra of $T$. Let $\d$ be a complement of $\n$ in $\t$, then the
complex vector space $\dc=\d+i\d$ is a complement of $\nc$ in $\tc$. We define a
projection $p \,\colon\, \dc\lorw \Tc/\Nc$ by the rule $p(Z)=\Pi(\exp Z),  Z\in \dc$,
where $\Pi \,\colon\, \Tc\lorw \Tc/\Nc$ denotes the canonical projection. Notice
that, by definition of $\exp$ the kernel of $p$ is a quasilattice $Q$ and $p$ induces
a group isomorphism $\dc/Q\cong \Tc/\Nc$. \qed We conclude this section with the
definition of holomorphic action of a complex quasitorus on a complex quasifold.
\begin{defn}[Holomorphic action]
    {\rm A {\em holomorphic action} of a complex quasitorus $\Dc$ on a complex quasifold $M$
    is a holomorphic mapping $\tau\,\colon\, \Dc\times M\longrightarrow M$
    such that $\tau(d_1\cdot d_2,m)=\tau(d_1,\tau(d_2,m))$ and
    $\tau(1_{\Dc},m)=m$ for all elements $d_1, d_2$ in the quasitorus $\Dc$ and for each point $m$ in
    the space $M$.}\end{defn}

\section{From simple polytopes to complex quasifolds}\label{tori}
Let $\d$ be a real vector space of dimension $n$, and let $\Delta$ be a simple convex
polytope of dimension $n$ in the dual space $\ddu$ (we recall that a polytope of
dimension $n$ is simple if there are exactly $n$ edges stemming from each of its
vertices). It is our intention to associate to this polytope a family of complex
quasifolds, in much the same way that one associates a toric variety to a simple
convex polytope that is rational. To do so we follow and extend the procedure for
constructing toric varieties as global quotients that is described by Guillemin in
\cite{vg}.

Write the polytope as
\begin{equation}\label{polydecomp}
\D=\bigcap_{j=1}^d\{\;\mu\in\ddu\;|\;\langle\mu,X_j\rangle\geq\lambda_j\;\}
\end{equation}
for some elements $\xd$ in the vector space $\d$ and some real numbers $\ld$. Let $Q$
be a quasilattice in the space $\d$ containing the elements $X_j$ (for example the
one that is generated by these elements) and let $\{\ed\}$ denote the standard basis
of $\rd$ and $\cd$; consider the surjective linear mapping
$$
\begin{array}{cccc}\label{pi}
\pi \,\colon\,& \R^d & \lorw & \d\\
    &   e_j& \longmapsto & X_j,
\end{array}
$$ and its complexification $$
\begin{array}{cccc}
\pic \,\colon\,& \C^d & \lorw & \dc\\
    &   e_j& \longmapsto & X_j.
\end{array}
$$ Consider the quasitorus $\d/Q$ and its complexification
$\dc/Q$. The mappings $\pi$ and $\pic$ each induce a group homomorphism, $$\Pi
\,\colon\, T^d=\rd/\Zd\lorw \d/Q$$ and $$\Pic \,\colon\, \Tdc=\cd/\Zd\lorw \dc/Q.$$
We define $N$ to be the kernel of the mapping $\Pi$ and $\Nc$ to be the kernel of the
mapping $\Pic$. Notice that neither $N$ nor $\Nc$ is a torus unless $Q$ is a honest
lattice. The mapping $\Pic$ defines an isomorphism
\begin{equation}\label{qtiso}
\Tdc/\Nc\longrightarrow \dc/Q
\end{equation}
\begin{remark}\label{polar}{\rm For the complexified group $\Nc$
the polar decomposition holds, namely
\begin{equation}\label{polareq}
\Nc=NA,
\end{equation}
where $A=\exp(i\n)$. In other words every element $w\in\Nc$ can be written {\em
uniquely} as $x\,\exp(i Y)$ where $x\in N$ and $Y\in\n$. This follows from the
definition of $N$ and $\Nc$, indeed
$\Nc=\{\,\exp(Z)\;|\;Z\in\C^d\quad\hbox{and}\quad\pic(Z)\in Q\,\}$. Write $Z=X+iY$,
then $\pic(Z)\in Q$ if and only if $\pi(X)\in Q$ and $\pi(Y)=0$, which implies
(\ref{polareq}).}
\end{remark}
Let $F$ denote a codimension-$k$ face of the polytope; this face is defined by a
system of $k$ equalities: $\langle\mu,X_j\rangle=\lambda_j$, for $j\in I\subset
\{1,\ldots,d\}$. Then we consider the $\Tdc$-orbit
$\cdf=\{\,\zd\in\cd\;|\;z_j=0\;\;\hbox{iff}\;\; j\in I\,\}$ and we take the union
over all the possible faces of the polytope
$$\cdd=\bigcup_F\cdf.$$ The group $\Nc$ acts on the space $\cdd$.
Consider the space of orbits $\cdd/\Nc$. We then have:
\begin{thm}\label{poltocx} Let $\d$ be a vector space of dimension $n$, and
let $\D\subset\ddu$ be a simple convex polytope. Choose inward-pointing normals to
the facets of $\D$, $X_1,\ldots,\X_d\in\d$, and let $Q$ be a quasilattice containing
these vectors. Then the corresponding quotient space $\cdd/\Nc$ is a complex
quasifold of dimension $n$. The complex quasitorus $\dc/Q$ acts on $\cdd/\Nc$, this
action is holomorphic and has a dense open orbit.
\end{thm}
Before we prove this theorem we need a lemma that will be crucial not only now but
also throughout the rest of the article. Consider any face $F$ of the polytope and
denote by $S^F_{\SC}$ the stabilizer of the $T^d_{\SC}$-action on $\C^d_F$ and by
$\s^F_{\SC}$ its Lie algebra. Then the fact that the polytope is simple implies that
\begin{equation} \nc\cap\s^F_{\SC}=\{0\},\label {free}\end{equation}
and this in turn implies that
\begin{equation} A\cap S^F_{\SC}=\{\mbox{Id}\},\label {Afree}\end{equation}
and that the following group is discrete
\begin{equation} \Nc\cap S^F_{\SC}=N\cap S^F_{\SC}.\label {Notfree}\end{equation}
Equation~(\ref{free}) also implies that $\pi$, restricted to $\s^F_{\SC}$, is
injective for any face $F$. In the special case that $F$ is equal to a vertex, $\mu$,
we get that $\smuc$ is a complement of $\nc$ in $\cd$, and that the linear mapping
$\pic$, when restricted to $\smuc$, defines a (very useful) isomorphism
\begin{equation}\label{keymap}\pi_{\mu}\,\colon\,\smuc\longrightarrow\dc.\end{equation}
One can then deduce:
\begin{lemma}\label{torustrick}Let $\mu$ be a vertex of the
polytope $\D$, consider the stabilizer $S^{\mu}_{\SC}$ of the orbit $\cd_{\mu}$, its
Lie algebra $\s^{\mu}_{\SC}$, and the discrete group $\G_{\mu}=S^{\mu}_{\SC}\cap\Nc$.
Then we have that\\ {\rm (i)} $T^d_{\SC}/\Smuc\simeq\Nc/\Gamma_{\mu}$;\\ {\rm (ii)}
$\Nc=\Gamma_{\mu}\exp{(\nc)}$;\\ {\rm (iii)} given any complement $\bc$ of $\smuc$ in
$\cd$, we have that
$$\nc=\{\,V-\pi_{\mu}^{-1}(\pic(V))\;|\;V\in\bc\,\}.$$
\end{lemma}\proof (i) Consider the group homomorphism $$
\begin{array}{cccc}
\lambda_{\mu} \,\colon\,& \Nc&\lorw&T^d_{\SC}/\Smuc\\
&n&\longmapsto&[n].
\end{array}
$$ Since $\nc$ and $\smuc$ are complementary, we have that
$T^d_{\SC}=\Smuc\Nc$, therefore $\lambda_{\mu}$ is surjective. The kernel of
$\lambda_{\mu}$ is given by $\G_{\mu}$, so $\lambda_{\mu}$ induces an isomorphism
$T^d_{\SC}/\Smuc\simeq\Nc/\Gamma_{\mu}$.\\(ii) Every element in $\Nc$ can be written
in the form $\exp{(Z)}$, where $Z\in\cd$ is such that $\pic(Z)\in Q$. Write now
$Z=Z-\pi_{\mu}^{-1}(\pic(Z))+\pi_{\mu}^{-1}(\pic(Z))$; it is easy to check that
$Z-\pi_{\mu}^{-1}(\pic(Z))\in\nc$, and that
$\exp{(\pi_{\mu}^{-1}(\pic(Z)))}\in\G_{\mu}$. The group $\G_{\mu}\cap\exp{(\nc)}$ is
not necessarily trivial, so the decomposition is not necessarily unique.
\\ (iii) Every element of the form $V-\pi_{\mu}^{-1}(\pic(V))$,
$V\in\bc$ clearly belongs to $\nc$. Conversely, write every element $Z\in\nc$ as
$Z=U+V$ according to the decomposition $\cd=\smuc\oplus\bc$, and notice that
$\pic(Z)=0$ implies that $U=-\pi_{\mu}^{-1}(\pic(V))$.\qed \vspace{.3cm}
\noindent{\mbox{\bf Proof of Theorem~\ref{poltocx}.\ \ }} We want to define a complex
quasifold atlas on the topological space $\cdd/\Nc$. We start by considering a rather
special covering of $\cdd$. To do so we restrict our attention to each vertex $\mu$ of
the polytope (the corresponding orbit $\cdmu$ has the smallest possible dimension
$d-n$) and we take the following $\Tdc$-invariant, open and connected neighborhood of
the orbit $\cdmu$ in $\cd$:
$$\vcmu=\{\,\zd\in\cd\;|\;z_j\neq0\;\;\hbox{if}\;\; j\notin
I\,\},$$ where $I$ is the index set corresponding to the vertex $\mu$ according to
the recipe given below Remark~\ref{polar}. Notice that we have $\cdd=\bigcup_{\mu}
\vcmu$, where $\mu$ ranges over all the vertices of the polytope $\Delta$. Indeed,
take $\zd\in\cdd$; then $\zd\in \C^d_G$ for some face $G$, and $\zd\in\vcmu$ for any
vertex $\mu$ contained in the face $G$. The opposite inclusion holds because the
polytope is simple. The neighborhoods $\vcmu$ are rather special, they are tubular.
Let us check this. Fix a point $\vz^{\mu}$ in the orbit $\cdmu$, for example
$\vz^{\mu}=(z^{\mu}_1,\cdots,z^{\mu}_d)$ with
\begin{equation}\label{zmu}\left\{
\begin{array}{ccc}
z^{\mu}_j=0 & \hbox{if} & j\in I\\ z^{\mu}_j=1 & \hbox{if} & j\notin I.
\end{array}
\right. \end{equation} Then there is a $\Tdc$-equivariant biholomorphism given by $$
\begin{array}{ccccc}
\Tdc&\times_{\Smuc}& \vtmu & \lorw & \vcmu\\
\left[t\right.&:&\left.\vz\right] & \longmapsto & t\cdot(\vz+\vz^{\mu}),
\end{array}
$$ where $$\vtmu=\prod_{i\in I}\C e_i\simeq \C^n$$ is the
holomorphic slice at $\vz^{\mu}$ for the $\Tdc$-action on $\cdd$ and
$$\Smuc=\{\,\zd\in \Tdc\;|\;z_j=1\;\;\hbox{if}\;\;j\notin I\,\}$$ is the stabilizer
of the action at $\vz^{\mu}$. Notice that $\vtmu=\smuc$. We now prove that the
subsets $\vmu=\vcmu/\Nc$ are complex charts for the quotient space $\cdd/\Nc$.
Consider the continuous surjective mapping $$
\begin{array}{cccc}
p_{\mu} \,\colon\,& \vtmu&\lorw&\vmu\\
&\vz&\longmapsto&[\vz^{\mu}+\vz].
\end{array}
$$ The discrete group $\Gamma_{\mu}=\Nc\cap\Smuc$ acts on the set
$\vtmu$ as follows: $$
\begin{array}{cccccc}
&\Gamma_{\mu}&\times& \vtmu&\longrightarrow&\vtmu\\ &(t &,&\vz)&\longmapsto& t\cdot
\vz.
\end{array}
$$ This action is holomorphic, free on a connected dense open set
and the mapping $p_{\mu}$ induces a continuous bijection
$$\phi_{\mu} \,\colon\,\vtmu/\Gamma_{\mu}\lorw \vmu.$$ To show
that $\phi_{\mu}$ is actually a homeomorphism, we only need to show that it is an
open mapping. This amounts to showing that, given a $\Gamma_{\mu}$-invariant open
subset $W$ of $\vtmu$, then $\Nc\cdot(\vz^{\mu}+W)$ is an open subset of $\cdd$.
Since $W$ is $\G_{\mu}$-invariant, by Lemma~\ref{torustrick}, (ii), we have that
$\Nc\cdot(\vz^{\mu}+W)=\exp(\nc)\cdot(\vz^{\mu}+W)$. Applying Lemma~\ref{torustrick},
(iii), we construct the surjective mapping
$$
\begin{array}{lcrcc}
W&\times&\Pi_{j\notin I} \C e_j  &\lorw&\exp(\nc)(\vz^{\mu}+W)\\
(\vw&,&V)&\longmapsto&\exp(V)\cdot\vz^{\mu}+\exp(-\pi_{\mu}^{-1}(\pic(V)))\cdot\vw.
\end{array}
$$ The determinant of its Jacobian matrix is given by
$\sum_{j=1}^{n}e^{K_j(Y)} +4\pi^2\sum_{h\notin I}e^{-4\pi Y_h}$, where
$Y=\frac{i}{2}(\overline{V}-V)$ and the $K_j$'s are real linear functions. This
implies, by the inverse function theorem, that $\exp(\nc)(\vz^{\mu}+W)$ is an open
subset of $\cdd$.

Let us now show that these charts are compatible. Consider two vertices of $\Delta$,
$\mu$ and $\nu$, and let $I$ and $J$ be the corresponding subsets of
$\{1,\ldots,d\}$. Assume that the corresponding charts, $\vmu$ and $\vnu$, have
non-empty intersection, and consider the two sets $\wmu=\fimu^{-1}(\vmu\cap\vnu)$ and
$\wnu=\finu^{-1}(\vmu\cap\vnu)$. We want to describe these two sets as submodels of
$\vtmu/\Gamma_{\mu}$ and $\vtnu/\Gamma_{\nu}$ respectively, and show that
$\gimunu=\finu^{-1}\circ\fimu$ is a biholomorphism of models. Consider
$\wtmu=p_{\mu}^{-1}(\fimu^{-1}(\vmu\cap\vnu))\subset\vtmu$ and
$\wtnu=p_{\nu}^{-1}(\finu^{-1}(\vmu\cap\vnu))\subset\vtnu$. Now, $$
\wtmu=\left(\prod_{j\in I\cap J}\C e_j \right)\times \left(\prod_{j\in I\setminus
I\cap J} \C^* e_j\right) $$ and
$$\wtnu=\left(\prod_{j\in I\cap J} \C e_j \right)\times
\left(\prod_{j\in J \setminus I\cap J} \C^*e_j\right). $$ Consider the universal
covers $\wfmu$, $\wfnu$ of $\wtmu$, $\wtnu$ respectively. Notice that we have that
$$\wfmu=\left(\prod_{j\in I\cap J}\C e_j \right)\times \left(\prod_{j\in I\setminus
I\cap J} \C e_j\right)\simeq\prod_{j\in I}\C e_j=\vtmu$$ and
$$\wfnu=\left(\prod_{j\in I\cap J}\C e_j= \right)\times
\left(\prod_{j\in J\setminus I\cap J} \C e_j\right)\simeq\prod_{j\in J}\C e_j=\vtnu,
$$ with projection maps $$
\begin{array}{ccc}
\wfmu&\lorw&\wtmu\\ (\vz,\vzeta)&\longmapsto&(\vz,\exp \vzeta)
\end{array}
$$ and $$
\begin{array}{ccc}
\wfnu&\lorw&\wtnu\\ (\vu,\veta)&\longmapsto&(\vu,\exp \veta).
\end{array}
$$ Consider the discrete groups $$\Lambda_{\mu}=\left\{\;(\exp
Z,W)\;|\;Z\in\bigoplus_{j\in I\cap J} \C e_j, \;W\in\bigoplus_{j\in I\setminus I\cap
J} \C e_j,\; \pic(Z+W)\in Q\;\right\} $$ and $$ \Lambda_{\nu}=\left\{\;(\exp
U,V)\;|\;U\in\bigoplus_{j\in I\cap J} \C e_j, \;V\in\bigoplus_{j\in J\setminus I\cap
J} \C e_j,\; \pic(U+V)\in Q\;\right\} $$ acting respectively on $\wfmu$, $\wfnu$ as
follows
$$
\begin{array}{lcccc}
(\Lambda_{\mu}&,&\wfmu)&\lorw&\wfmu\\ ((\exp Z,W)&,&(\vz,\vzeta))&\longmapsto&(\exp
Z\cdot \vz, W + \vzeta)
\end{array}
$$ $$
\begin{array}{lcccc}
(\Lambda_{\nu}&,&\wfnu)&\lorw&\wfnu\\ ((\exp U,V)&,&(\vu,\veta))&\longmapsto&(\exp
U\cdot \vu, V + \veta).
\end{array}
$$ Notice that the projections induce homeomorphisms
$\wfmu/\Lambda_{\mu}\simeq\wmu$ and $\wfnu/\Lambda_{\nu}\simeq\wnu$. Now we want to
show that there is an equivariant biholomorphism $\gimunus\colon\wfmu\lorw \wfnu$
that projects down to the mapping $\gimunu$. Consider the isomorphisms $\pi_{\mu}$
and $\pi_{\nu}$ defined by (\ref{keymap}). Notice that $\pi_{\nu}^{-1}\cdot\pi_{\mu}$
defines a biholomorphism from $\wfmu$ to $\wfnu$ that is equal to the identity on
$\prod_{j\in I\cap J}\C e_j$. Let moreover $\rho$, respectively $\sigma$, denote the
projection of $\wfnu$ onto the factor $\prod_{j\in I\cap J}\C e_j$, respectively
$\prod_{j\in J\setminus I\cap J}\C e_j$. Define $$
\begin{array}{cccc}
\gimunus\colon&\wfmu&\lorw&\wfnu\\
&(\vz,\vzeta)&\longmapsto&\left(\exp
\left(\rho\left(\pi_{\nu}^{-1}\circ\pi_{\mu}(\vzeta)\right)\right)\cdot
\vz,\sigma\left(\pi_{\nu}^{-1}\circ\pi_{\mu}(\vzeta)\right)\right).
\end{array}
$$ It is straightforward to check that $\gimunus$ is an
equivariant biholomorphism that projects down to $\gimunu$. Now complete this
collection of charts with all the other compatible charts. This concludes the proof
that $\cdd/\Nc$ is a complex quasifold of dimension $n$. The standard holomorphic
action of $\Tdc$ on $\cdd$ induces a holomorphic action of $\Tdc/\Nc$ on $\cdd/\Nc$.
The isomorphism (\ref{qtiso}) allows us to define a holomorphic action
$\tau\,\colon\,(\dc/Q)\times(\cdd/\Nc)\rightarrow \cdd/\Nc$. To check that $\tau$ is
holomorphic, consider, for each vertex $\mu$ of the polytope, the chart $\vmu$, and
the bijection $\pi_{\mu}$ given in (\ref{keymap}). We then have that the following
diagram commutes $$
\begin{array}{ccc}
\dc\times \vtmu & \stackrel{\tilde{\tau}_{\mu}}{\mbox{\LARGE $\longrightarrow$}} &
\vtmu \\ (W,\vz)&{\mbox{\LARGE
$\longmapsto$}} & \exp (\pi_{\mu}^{-1}(W))\cdot \vz\\
{}\!\!\!\!\!\stackrel{}{}\,\stackrel{}{\mbox{\LARGE $\downarrow$}}
&  & \stackrel{}{\mbox{\LARGE $\downarrow$}}\,\stackrel{}{}\\
(\dc/Q)\times(\vtmu/\gmu) & \stackrel{\tau_{\mu}}{\mbox{\LARGE $\longrightarrow$}} &
\vtmu/\gmu\\ ([W],[\vz])&{\mbox{\LARGE $\longmapsto$}} & [\exp
(\pi_{\mu}^{-1}(W))\cdot \vz]\\{}\!\!\!\!\!\stackrel{}{}\,\stackrel{}{\mbox{\LARGE
$\downarrow$}} &  & \stackrel{}{\mbox{\LARGE $\downarrow$}}\,\stackrel{}{}\\
(\dc/Q)\times\vmu & \stackrel{\tau}{\mbox{\LARGE $\longrightarrow$}} &
\vmu\\([W],[\vz+\vz_\mu])&{\mbox{\LARGE $\longmapsto$}} & [\exp
(\pi_{\mu}^{-1}(W))\cdot ( \vz +
\vz_\mu)]=[\exp(\pi_{\mu}^{-1}(W))\cdot\vz+\vz_\mu]\\
\end{array}
$$ and that $\tilde{\tau}_\mu$ is a holomorphic mapping. Notice
finally that the dense open orbit for this action is given by $\cd_F/\Nc$, where
$F=\mbox{Int}(\D)$.\qed
\begin{remark}{\rm Suppose that the polytope
$\D$ is rational with respect to a lattice $L\subset\d$. Choose inward-pointing
normal vectors to the facets of $\D$ that are primitive in $L$, and take $Q=L$. Then
the corresponding quotient $\cdd/\Nc$ is a complex orbifold and contains a dense open
orbit for the action of the honest torus $\Tc=\dc/L$; it is the usual toric variety
associated to the polytope $\D$.}
\end{remark}
\begin{remark}{\rm Notice that, just like for honest toric varieties,
the quasitorus $\dc/Q$ is contained in $\cdd/\Nc$ as a dense open subset, and the
action of $\dc/Q$ on $\cdd/\Nc$ extends the standard action of $\dc/Q$ on
itself.}\end{remark} Let us test this construction on three examples. Only the third
is an example of a truly non-rational polytope.
\begin{ex}[The unit interval]\label{quasisfera}{\rm
Let us consider the polytope $[0,1]\subset\rdu$ and let us apply the above
construction to the choice of vectors $X_1=s$, $X_2=-t$, with $s, t >0$,
$s/t\notin\Q$, and to the choice of quasilattice $Q=s\Z+t\Z$. It is easy to check
that $\cddtwo=\ctwo \setminus\{0\}$ and that $\Nc= \{\,(e^{2\pi iZ},e^{2\pi
i\frac{s}{t}Z})\;|\; Z\in\C\,\}$. We cover $\cddtwo/\Nc$ with two charts, one for
each of its vertices, $\sc{s}=\{0\}$ and $\sc{n}=\{1\}$: $V_{\sc s}
=\left\{\,[z_1:z_2]\in\cddtwo/\Nc\;|\; z_2\neq 0\,\right\}$ and $V_{\sc
n}=\left\{\,[z_1:z_2]\in\cddtwo/\Nc\;|\; z_1\neq 0\,\right\}$. Consider now the the
discrete group $\G_{\sc s}=\Z$ acting on the set $\vts=\C$ according to the rule
$(k,z)\mapsto e^{2\pi i \frac{t}{s}k} \cdot z$; this action is holomorphic and free
on the connected, dense subset $\vts\setminus \{0\}$ and the mapping
$$\begin{array}{cccc}\phi_{\sc s}\,\colon\,&\vts/\G_{\sc
s}&\longrightarrow &V_{\sc s}\\ &[z]&\longmapsto &[z:1]\end{array}$$ is a
homeomorphism. Similarly the group $\G_{\sc n}=\Z$ acts on $\vtn=\C$ by the rule
$(m,z)\mapsto e^{2\pi i \frac{s}{t}m} \cdot w$; this action is holomorphic and free
on the connected, dense subset $\vtn\setminus\{0\}$ and the mapping
$$\begin{array}{cccc}\phi_{\sc n}\,\colon\,&\vtn/\G_{\sc n}&\longrightarrow &V_{\sc
n}\\ &[w]&\longmapsto &[1:w]\end{array}$$ is a homeomorphism. Let us show that these
two charts are compatible. The set $\phi_{\sc s}^{-1}(\vs\cap\vn)$ is a submodel of
$\vts/\G_{\sc s}$: it is the quotient of the space $\wts=\C$ by the action of
$\Lambda_{\sc s}=\Z^2$ given by $(h,k)\cdot \zeta=\zeta+h+ \frac{t}{s}k$. Similarly
the set $\phi_{\sc n}^{-1}(\vs\cap\vn)$ is the quotient of the space $\wtn=\C$ by the
action of $\Lambda_{\sc n}=\Z^2$ given by $(l,m)\cdot \eta=\eta+l+\frac{s}{t}m$. The
bijective mapping
\begin{eqnarray*}
g_{\sc sn}=\phi_{\sc n}^{-1}\circ\phi_{\sc s}\,\colon \phi_{\sc s}^{-1}(V_{\sc s}\cap
V_{\sc n}) &\longrightarrow &\phi_{\sc n}^{-1}(V_{\sc s}\cap V_{\sc n})\\
\left[z\right] &\longmapsto &\left[w=z^{-\frac{s}{t}}\right]
\end{eqnarray*}
is a biholomorphism of models: its lift is given by $\zeta\longmapsto
\eta=-\frac{s}{t}\zeta$. Now complete this collection with all the other compatible
charts.}
\end{ex}
\begin{remark}{\rm
We can see from the previous example that complex quasifolds corresponding to the
same polytope $\Delta$ are not in general biholomorphic, in fact they are not even
diffeomorphic. This was also visible in the rational case: the same construction
applied to $s$, $t$ relatively prime integers yields in fact either ordinary or
weighted projective space,  which have different complex orbifold structures. This
may appear to be in contradiction with Theorem 9.4 of Lerman-Tolman \cite{lt}, which
implies that these spaces are biholomorphic - in reality their notion of
biholomorphism is algebraic and, unlike ours, does not keep track of the orbifold
structure.}
\end{remark}
\begin{ex}[The right triangle]\label{triangle}{\rm
Consider the right triangle in $\rtwodu$ of vertices $(0,0)$, $(s,0)$ and $(0,t)$,
where $s, t$ are two positive real numbers such that $s/t\notin\Q$. We apply the
construction to the choice of inward-pointing normals $X_1=(1,0)$, $X_2=(0,1)$,
$X_3=(-t,-s)$ and to the quasilattice $Q=X_1\Z+X_2\Z+X_3\Z$. Then we have that
$\cddth=\cthree\setminus\{0\}$ and $\Nc=\{\,(e^{2\pi i t Z}, e^{2\pi i s Z}, e^{2\pi
i Z})\;|\; Z\in\C\,\}$. We cover $\cddth/\Nc$ with three charts, one for each of its
vertices: $V_1 =\left\{\,[z_1:z_2:z_3]\in\cddth/\Nc\;|\; z_1\neq 0\,\right\}$,
$V_2=\left\{\,[z_1:z_2:z_3]\in\cddth/\Nc\;|\; z_2\neq 0\,\right\}$ and
$V_3=\left\{\,[z_1:z_2:z_3]\in\cddth/\Nc\;|\; z_3\neq 0\,\right\}$. The discrete
group $\G_1=\Z$ acts on $\vtone=\ctwo$ according to the rule
$h\cdot(z_2,z_3)=(e^{2\pi i \frac{s}{t} h}\cdot z_2, e^{2\pi i \frac{1}{t}h}\cdot
z_3)$ and the mapping
$$\begin{array}{cccc}\phi_1\,\colon\,&\vtone/\G_1&\longrightarrow
&V_1\\ &[z_2,z_3]&\longmapsto &[1:z_2:z_3]\end{array}$$ is a homeomorphism. Similarly
the group $\G_2=\Z$ acts on $\vttwo=\ctwo$ by the rule $k\cdot(z_1,z_3)=(e^{2\pi i
\frac{t}{s} k}\cdot z_1, e^{2\pi i \frac{1}{s}k}\cdot z_3)$ and the mapping
$$\begin{array}{cccc}\phi_2\,\colon\,&\vttwo/\G_2&\longrightarrow
&V_2\\ &[z_1,z_3]&\longmapsto &[z_1:1:z_3]\end{array}$$ is a homeomorphism. Finally
the group $\G_3=\Z$ acts on $\vtthree=\ctwo$ according to the rule $l\cdot (z_1,z_2) =
(e^{2\pi i t l}\cdot z_1,e^{2\pi i s l}\cdot z_2)$ and the mapping
$$\begin{array}{cccc}\phi_3\,\colon\,&\vtthree/\G_3&\longrightarrow
&V_3\\ &[z_1,z_2]&\longmapsto &[z_1:z_2:1]\end{array}$$ is a homeomorphism. The
changes of charts work as in the previous example.}
\end{ex}
\begin{ex}[The regular pentagon]{\rm
Let us take the regular pentagon in $\rtwodu$ of vertices $(1,0,)$, $(a,b)$, $(c,d)$,
$(c,-d)$ and $(a,-b)$, where $a=\cos{\frac{2\pi}{5}}$, $b=\sin{\frac{2\pi}{5}}$,
$c=\cos{\frac{4\pi}{5}}$, $d=\sin{\frac{4\pi}{5}}$. There exists no lattice $L$ with
respect to which this simple convex polytope is rational. Let us apply the above
construction to the choice of inward-pointing normal vectors $X_1=(1,0)$, $X_2=(a,b)$,
$X_3=(c,d)$, $X_4=(c,-d)$, $X_5=(a,-b)$ and to the choice of quasilattice
$Q=\sum_{j=1}^{5}X_j\Z$. Then we have that
$$\Nc=\left\{\,\left(e^{2\pi i Z_1},e^{2\pi i Z_2},e^{2\pi i
Z_3},e^{2\pi i[2a(Z_2-Z_3)+Z_1)]}, e^{2\pi i[2a(Z_2-Z_1)+Z_3)]}\right)\;|\;
(Z_1,Z_2,Z_3) \in\cthree \,\right\}.$$ We cover the quasifold $\cddfv/\Nc$ with five
charts: $$V_1=\left\{\,[z_1:z_2:z_3:z_4:z_5]\in\cddfv/\Nc\;|\; z_1, z_2, z_3 \neq
0\,\right\},$$
$$V_2=\left\{\,[z_1:z_2:z_3:z_4:z_5]\in\cddfv/\Nc\;|\; z_2, z_3,
z_4 \neq 0\,\right\},$$
$$V_3=\left\{\,[z_1:z_2:z_3:z_4:z_5]\in\cddfv/\Nc\;|\; z_3, z_4,
z_5 \neq 0\,\right\},$$
$$V_4=\left\{\,[z_1:z_2:z_3:z_4:z_5]\in\cddfv/\Nc\;|\; z_1, z_4,
z_5\neq 0\,\right\}$$ and $$V_5 =\left\{\,[z_1:z_2:z_3:z_4:z_5]\in\cddfv/\Nc\;|\;
z_1, z_2, z_5\neq 0\,\right\}.$$ The discrete group $\G_1=\Z^3$ acts on the set
$\vt_1=\ctwo$ according to the rule $(h,k,l)\cdot(z_4,z_5)=(e^{4\pi i a(k-l)}\cdot
z_4, e^{4\pi i a(k-h)}\cdot z_5)$ and the mapping
$$\begin{array}{cccc}\phi_1\,\colon\,&\vt_1/\G_1&\longrightarrow
&V_1\\ &[z_4:z_5]&\longmapsto &[1:1:1:z_4:z_5]\end{array}$$ is a homeomorphism. We
leave the discussion of the other four charts, and their mutual compatibility, to the
reader.}
\end{ex}

\section{K\"ahler structures}
In the previous section we associated to any simple convex polytope
$\Delta\subset\ddu$, together with a choice of normal vectors $X_j$, $j=1,\ldots,d$,
and of a quasilattice $Q\subset\d$ containing the $X_j$'s, a complex quasifold of
dimension $n$. Starting with the same data it is possible to construct a symplectic
quasifold of dimension $2n$, endowed with an effective Hamiltonian action of the
quasitorus $\d/Q$, such that the image of the corresponding moment mapping is exactly
$\Delta$ (see \cite[Theorem 3.3]{p}). The complex and symplectic quasifolds are both
described as orbit spaces, the first is the quotient of $\cdd$ by $\Nc$, the second
is a symplectic quotient with respect to the action of $N$ on $\C^d$. We want to
prove that the complex and symplectic quotient can be identified, according to a
general principle initiated by Kempf-Ness \cite{kempfness} and later developed by
Kirwan \cite{kirwan} and Ness \cite{ness}. More precisely, we shall prove that the
two quotients are diffeomorphic and that the complex and symplectic structure are
compatible, and thus define the structure of a K\"ahler quasifold (see
Definition~\ref{kahform}).

To begin with, let us briefly recall from \cite{p} the construction of the symplectic
quotient that we are interested in. Consider the mapping $J(\vz)=\sum_{j=1}^d
(|z_j|^2+\lambda_j)e_j^*$, where the $\lambda_j$'s are given in (\ref{polydecomp})
and are uniquely determined by our choice of normal vectors. The mapping $J$ is a
moment mapping for the standard action of $T^d$ on $\C^d$. Consider now the subgroup
$N\subset T^d$ and the corresponding inclusion of Lie algebras
$\iota\,\colon\,\n\rightarrow\R^d$. The mapping $\Psi\,\colon\,\C^d\rightarrow \n^*$
given by $\Psi={\iota}^*\circ J$ is a moment mapping for the induced action of $N$ on
$\cd$. Then the quotient space $\squo$ is a compact symplectic quasifold of dimension
$2n$; the quasitorus $\d/Q$ acts on $\squo$ in an effective, Hamiltonian fashion, and
the image of the corresponding moment map is the polytope $\D$. Let us now define a
mapping between the symplectic and complex quotient. To define this mapping and to
show that it is bijective we adapt to our setting the method described by Guillemin
in \cite[Appendix 1]{vg} for the smooth case. Before we can define this mapping, we
need a preliminary lemma.
\begin{lemma} The zero set $\Psi^{-1}(0)$ is contained in $\cdd$.
Moreover, for any face $F$ of the polytope $\Delta$, the orbit $\cdf$ intersects
$\Psi^{-1}(0)$ in at least one point. \label{contenuto}
\end{lemma}
\proof We summarize the proof given in \cite[page 115]{vg}, which goes through
without modification. Consider the exact sequence
\begin{equation}
\label{exactsequence}
0\lorw\d^*\stackrel{\pi^*}{\lorw}(\R^d)^*\stackrel{{\iota}^*}{\lorw}(\n)^*\lorw0,
\end{equation}
where $\pi$ is the projection defined in (\ref{pi}). Notice that $\vz\in\zset$ if,
and only if ${\iota}^*(J(\vz))=0$. By (\ref{exactsequence}) we have that, for a given
$\vz\in\zset$, there exists a unique $\zeta\in\d^*$ such that $J(\vz)=\pi^*(\zeta)$.
By making use of the explicit expression of $J$ we find
\begin{equation}
\label{formula} |z_j|^2=\langle \zeta,X_j\rangle-\lambda_j,\quad\quad j=1,\ldots, d.
\end{equation}
This implies that, given $\vz\in\zset$, the corresponding $\zeta$ lies in $\D$ and
that, again by (\ref{formula}), $\vz$ is in $\C^d_F$, where $F$ is the face of
$\Delta$ containing $\zeta$ in its interior. Similarly, given any face $F$ of the
polytope, we can always find a point $\vz\in\zset\cap\C^d_F$. \qed By Lemma
\ref{contenuto}, there is an injection
\begin{equation}
{\cal I}\,\colon\,\Psi^{-1}(0)\longrightarrow\cdd,
\end{equation}
which induces the mapping
\begin{equation}
\label{lachi} \chi\,\colon\, \squo \lorw \cquo
\end{equation}
that sends an $N$-orbit to the corresponding $\Nc$-orbit. This mapping is equivariant
with respect to the actions of the quasitori $\d/Q$ and $\dc/Q$. We are now ready to
state the main result of this section.
\begin{thm}
\label{teoremadellachi} Let $\d$ be a vector space of dimension $n$, and let
$\D\subset\ddu$ be a simple convex polytope. Choose inward-pointing normals to the
facets of $\D$, $X_1,\ldots,\X_d\in\d$, and let $Q$ be a quasilattice containing
these vectors. Then the mapping $$\chi\,\colon\,\squo\lorw\cquo$$ is an equivariant
diffeomorphism of quasifolds. Moreover the induced symplectic form on the complex
quasifold $\cquo$ is K\"ahler.
\end{thm}
For the definitions of symplectic quasifold, diffeomorphism of quasifolds, and
pullback of a differential form, we refer the reader to \cite{p}. Before we can
proceed with the proof of this theorem we need a number of remarks and lemmas.
\begin{remark}\label{movesout}{\rm The group $A=\exp(i\n)$ moves any point $\vz\in\zset$ out of
$\zset$. To see this, remark that, for any vector $Y\in\n$, the corresponding
function $\Psi_{Y}$ satisfies $d\Psi_{Y}={\imath}_{\tilde Y}\omega_0$, where
$\omega_0=\frac{1}{2\pi i}\sum_{j=1}^d dz_j\wedge d\overline{z}_j$ is the standard
symplectic form of $\C^d$ and $\tilde{Y}$ is the vector field generated by $Y$. This
implies that the gradient of the function $\Psi_{Y}$ is the vector field generated by
$iY$, which by (\ref{Afree}) does not vanish in $\cdd$. The function $\Psi_{Y}$ is
therefore strictly increasing along its gradient flow, hence
$\exp(tiY)\cdot\vz\notin\zset$, for all $Y\in\n$ and all $t\neq 0$.}
\end{remark}
\begin{lemma}
The following facts are equivalent:\\{\rm (i)} the mapping $\chi$ is bijective;\\{\rm
(ii)} every $\Nc$-orbit in $\cdd$ intersects $\Psi^{-1}(0)$ in an $N$-orbit;\\{\rm
(iii)} every $A$-orbit in $\C^d_{\Delta}$ intersects $\Psi^{-1}(0)$ in at least one
point. \label{equivalenze}
\end{lemma}
\proof It is obvious that (i) is equivalent to (ii). Before going
on with the proof let us re-state the second and third point:\\
(ii) for each $\Nc$-orbit through a point $\vz\in\cdd$ there exists a point
$\vu\in\Psi^{-1}(0)$ such that $(\Nc\cdot \vz)\cap\Psi^{-1}(0)=N\cdot\vu$;\\ (iii)
for each $A$-orbit through a point $\vz\in\cdd$ there exists at least one point
$\vw\in (A\cdot\vz)\cap\Psi^{-1}(0)$.

The proof is based on three key facts: the polar decomposition for the group $\Nc$
(see Remark~\ref{polar}), the invariance of $\Psi^{-1}(0)$ under the action of $N$,
and Remark~\ref{movesout}. To see that (ii) implies (iii), consider an $A$-orbit
through $\vz\in\cdd$. Then, by (ii), there exists $\vu\in\Psi^{-1}(0)$ such that
$(\Nc\cdot\vz)\cap\Psi^{-1}(0)=N\cdot\vu$. Since $\Nc=NA$, there exists $\vw\in
N\cdot\vu\subset\Psi^{-1}(0)$ such that $\vw\in A\cdot\vz$. To see that (iii) implies
(ii), consider any $\Nc$-orbit through $\vz$. By assumption, the orbit $A\cdot\vz$
intersects $\Psi^{-1}(0)$ in at least one point $\vw$. Since $\Nc=NA$ we have that
$\Nc\cdot\vz=\Nc\cdot\vw$, that $N\cdot\vw\subset(\Nc\cdot\vz)\cap\Psi^{-1}(0)$, and
that this inclusion is in fact an equality by Remark~\ref{movesout}.\qed
\begin{remark}\label{unicita}{\rm
If there exists a point $\vw\in (A\cdot\vz)\cap\zset$, then $\vw$ is uniquely
determined by this property. If $\vw$ and $\vw'$ are in $(A\cdot\vz)\cap
\Psi^{-1}(0)$, then there exists an $a\in A$ such that $a\cdot\vw=\vw'$. By
Remark~\ref{movesout} we have that $a$ is the identity and that $\vw=\vw'$.}
\end{remark}
The following result will be essential to proving that the mapping $\chi$ is a
bijection.
\begin{lemma}\label{legendre}
The moment mapping $\Psi$ maps any $A$-orbit diffeomorphically onto an open convex
cone in $\n^*$; moreover, if two $A$-orbits lie in the same $T^d_{\SC}$-orbit
$\cdf\subset\cdd$, then their images with respect to $\Psi$ are identical.
\end{lemma}
\proof The proof given in \cite[Appendix~1: Theorem~2.1 and Theorem~2.2]{vg} applies,
the only delicate point here is to notice that, by (\ref{Afree}), $A$ still acts
freely on $\cdd$. Notice also that, by (\ref{Notfree}), this is certainly not true
for the action of $N$. We now want to study the image by $\Psi$ of an $A$-orbit
through a point $\vz\in\cdd$, namely
$$\Psi(A\cdot\vz)=\left\{\,\Psi\left(e^{-2\pi \alpha_1(X)}\cdot
z_1,\ldots, e^{-2\pi \alpha_d(X)}\cdot z_d\right)\;|\;X\in\n\,\right\},$$ where
$\alpha_j=\iota^*(e_j^*)$, $j=1,\ldots,d$. Since $A$ acts freely on $\cdd$, the
exponential mapping defines a diffeomorphism between $\n$ and $A\cdot z$, and we can
identify the $A$-orbit with $\n$. Therefore the set we are interested in is the image
of
the mapping $$\begin{array}{cccc} f \,\colon\,&\n&\lorw&\n^*\\
&X&\longmapsto&\Psi\left(e^{-2\pi \alpha_1(X)}\cdot z_1,\ldots,e^{-2\pi
\alpha_d(X)}\cdot z_d\right).
\end{array}
$$ The point $\vz$ lies in $\C^d_F$ for a face $F$ of the polytope
$\D$. Let $I$ be the corresponding set of indices, as defined in Section~2, then
$$f(X)=\sum_{j\notin I}e^{-4\pi\alpha_j(X)}|z_j|^2\alpha_j+\lambda_j\alpha_j.$$ Now we
only have to check that
\begin{equation}
\label{spannstar} \mbox{span}\{\alpha_j\in\n^*\;|\;j\notin I\}=\n^*;
\end{equation}
this is equivalent to saying that there does not exist a non zero $X\in\n$ such that
$\alpha_j(X)=0$ for all $j\notin I$. Indeed, if there was such an $X$, it would lie
in $\n\cap\s^F_{\SC}$ which is $\{0\}$. Using (\ref{spannstar}) one can prove,
following \cite{vg}, that $f$ maps $\n$ diffeomorphically into an open convex cone of
$\n^*$, and also that the image of $f$ depends only on $I$. \qed We are now ready to
proceed with the proof of the
main result of this section.\\

\noindent{\mbox{\bf Proof of Theorem~\ref{teoremadellachi}.\ \ }}The proof is divided
into several steps. We first prove that the mapping $\chi$ is bijective and
continuous. Then we lift $\chi$ locally to prove that it is a diffeomorphism, and
finally we show that the pull-back via $\chi^{-1}$ of the symplectic form of $\squo$,
is K\"ahler on the complex quasifold $\cquo$.

In order to prove that $\chi$ is bijective we will show that Lemma~\ref{legendre}
implies that, for a given $A$-orbit through a point $\vz\in\cdd$, there exists a
point $\vw\in\Psi^{-1}(0)$ such that $(A\cdot\vz)\cap\Psi^{-1}(0)=\{\vw\}$. Then we
will apply Lemma~\ref{equivalenze} to conclude. To prove the existence of $\vw$,
consider the $T^d_{\SC}$-orbit $\C^d_F$ containing $A\cdot\vz$. By Lemma
\ref{contenuto} this $\C^d_F$ intersects $\zset$ in at least one point, $\vu$. On the
other hand $\Psi(A\cdot\vz)=\Psi(A\cdot\vu)\subset\n^*$  by Lemma~\ref{legendre}, and
the common image contains $0=\Psi(\vu)$. So there exists $a\in A$ such that
$\Psi(a\cdot \vz)=0$; then $\vw=a\cdot \vz\in (A\cdot\vz)\cap\Psi^{-1}(0)$. Notice
that, by Remark~\ref{unicita}, the point $\vw$ is uniquely determined. Continuity of
the mapping $\chi$ is implied by that of the immersion $\cal I$.

For the next step we need a covering of the symplectic quotient by a collection of
its charts. We refer the reader to \cite[Theorem~3.1]{p} for the construction of the
full atlas. Let $\mu$ be a vertex of the polytope $\D$, and let
$I=\{r_1,\ldots,r_n\}\subset \{1,\ldots,d\}$ denote the corresponding subset of
indices. For each vertex of $\D$ we want to construct a corresponding chart. Let
$(a_{jh})\in M_{(n,d)}(\R)$ denote the matrix of the mapping
$\pi\,\colon\,\rd\rightarrow\d$ with respect to the basis $\{\,X_j\,|\, j\in I\,\}$
of $\d$ and the standard basis of $\rd$. Consider the open subset of $\vtmu$ defined
as follows
$$\utmu=\left\{ \vw\in\vtmu \,\left|\right.\, \sum_{j=1}^{n}
a_{jh} \left(|w_{r_j}|^2+\lambda_{r_j}\right)-\lambda_h>0, \, h\notin I\,\right\}.$$
The group $\G_{\mu}$ defined in Lemma~\ref{torustrick} acts on $\utmu$. Consider, for
a given element $\vw\in \utmu$, the element
$\vw^{\mu}=(w^{\mu}_1,\cdots,w^{\mu}_d)\in\cd$ defined as follows
$$ \left\{
\begin{array}{lcc}
w^{\mu}_h=0 & \hbox{if} & h\in I\\ w^{\mu}_h=\sqrt{\sum_{j=1}^{n} a_{jh}
\left(|w_{r_j}|^2+\lambda_{r_j}\right)-\lambda_h} & \hbox{if} & h\notin I.
\end{array}
\right.$$ Notice that $\vw+\vw^{\mu}$ belongs to $\zset$: define $\nu\in\d^*$ so that
$<\nu,X_j>=|w_j|^2+\lambda_j$, for all $j\in I$; then it is easy to check that
$J(\vw+\vw^{\mu})=\pi^*(\nu)$ and $\nu\in\Delta$. Consider now the open sets
$\ucmu=\vcmu\cap\zset\subset\zset$ and $\umu=\ucmu/N\subset\zset/N$. Then the
surjective mapping $$
\begin{array}{cccc}
q_{\mu} \,\colon\,& \utmu&\lorw&\umu\\
&\vw&\longmapsto&[\vw^{\mu}+\vw]
\end{array}
$$ induces a homeomorphism $$
\begin{array}{cccc}
\psi_{\mu} \,\colon\,& \utmu/\G_{\mu}&\lorw&\umu\\
&[\vw]&\longmapsto&[\vw^{\mu}+\vw].
\end{array}
$$ The above data define a chart, and the union of the $\umu$'s,
for $\mu$ ranging over all the vertices of $\D$, cover the symplectic quotient.

Let us now show that the mapping $\chi$ lifts to a diffeomorphism
$\tilde{\chi}_{\mu}\,\colon\,\utmu\rightarrow\vtmu$ for each vertex $\mu$. For each
$\vw\in\utmu$, there exists a unique element $a(\vw)\in A$ such that $a(\vw)\cdot
(\vw+\vw^{\mu})$ is of the form $\vz+\vz^{\mu}$, where $\vz\in\vtmu$, and $\vz^{\mu}$
is given by (\ref{zmu}). Then we define $\tilde{\chi}_{\mu}(\vw)=\vz$. We compute
$\vz$ explicitly following Lemma~\ref{torustrick}:
$\vz=\tilde{\chi}_{\mu}(\vw)=\exp{\left(-i\pi_{\mu}^{-1}(\pic(C(\vw))\right)}
\cdot\vw$, where the $h$-component of $C(w)\in\prod_{h\notin J} \C e_h$ is given by
$\frac{1}{4\pi}\log{A_h}$, with $A_h=\sum_{j=1}^{n}
a_{jh}(|w_{r_j}|^2+\lambda_{r_j})-\lambda_h$. Notice that the mapping $\chitmu$ thus
defined is equivariant with respect to the action of $\G_{\mu}$ on $\utmu$ and
$\vtmu$, and it is a lift of $\chi$. It is now straightforward to deduce the
bijectivity of $\chitmu$ from that of $\chi$, using the explicit expression for
$\chitmu$. In order to prove that $\chitmu$ is a diffeomorphism we apply the inverse
function theorem. Therefore we only have to check that its Jacobian matrix
$D(\chitmu)$ is non-degenerate on $\utmu$. Let $(\underline{x},\underline{y})$ be
real coordinates in $\utmu$, it turns out that the Jacobian matrix has the following
form $$ D(\chitmu)(\vx,\vy)=e^r\left(I_{2n}+\left(\begin{array}{c}M\\N
\end{array}\right)(^t\!M,^t\!N)
\right),$$ where $r$ is the function $\sum_{h\notin I }a_{kh}\frac{1}{2}\log{A_h}$,
and where $M,N$ are $(n,d-n)$ matrices with entries $u_{k}a_{kh}\frac{1}{\sqrt{A_h}}$
and $v_{k}a_{kh}\frac{1}{\sqrt{A_h}}$ respectively. This implies that $D(\chitmu)(w)$
is positive definite for every $\vw\in\utmu$. To conclude that $\chi$ is a
diffeomorphism observe that the continuity of the equivariant mapping $\chitmu^{-1}$
for each vertex $\mu$ implies that $\chi^{-1}$ is continuous, since $\phi_{\mu}$ and
$\psi_{\mu}$ are homeomorphisms.

Now, having proved that $\chi$ is a diffeomorphism, we can consider the complex
quotient endowed with a symplectic and a complex structure. We want to prove that the
symplectic form is K\"ahler. This can be checked pointwise: let
$\hat{\vw}=\vw+\vw^{\mu}$ be a point in $\ucmu$. Then the $N$-orbit through
$\hat{\vw}$ is contained in $\zset$ and the $A$-orbit through $\hat{\vw}$ is
orthogonal to $\zset$. This gives the isomorphism
$$T_{\hat{\vw}}(\zset)/T_{\hat{\vw}}(N\cdot\hat{\vw})\simeq
\C^d/T_{\hat{\vw}}(\Nc\cdot\hat{\vw}).$$ To conclude the proof we only have to remark
that the symplectic structure on $\squo$ and the complex structure on $\cquo$, read on
$T_{\hat{\vw}}(\zset)/T_{\hat{\vw}}(N\cdot\hat{\vw})$, are exactly the ones induced
by the standard complex and symplectic structures of $\C^d$. \qed
\begin{remark}\label{manykhaler}{\rm
Consider the complex quasifold $\cdd/\Nc$ constructed in Theorem \ref{poltocx}.
Notice that, by varying the coefficients $\lambda_j$ in (\ref{polydecomp}), we can
produce many simple convex polytopes, each allowing the same choice of inward-pointing
normals $X_j$ and of quasilattice $Q$ that we had made for $\D$ (for example we can
``inflate'' $\Delta$). The corresponding complex quasifolds are therefore exactly the
same, but on the symplectic side in general one obtains non-equivalent symplectic
structures. In other terms we obtain many non-isometric K\"ahler structures on the
same complex quasifold $\cdd/\Nc$.}
\end{remark}
We conclude this section with the discussion of an example. The unit interval
provides a very useful model for a thorough understanding of the diffeomorphism
$\chi$ and a neat example of K\"ahler quasifold.
\begin{ex}\label{quasisferaii}
{\rm Consider the polytope $[0,1]\subset\rdu$, with the same choice of vectors and
quasilattice made in Example \ref{quasisfera} (we refer to that example for the
notation). We can associate to these data the complex quasifold $\ctwo_{\D}/\Nc$ and
also a symplectic quasifold, the {\em quasisphere} of \cite[Examples~1.13, 1.19, 2.10
and 3.5]{p}, of which we recall only some features.

Consider the symplectic quasifold $\squo$, where
$$\zset=\{\,(z,w)\in\ctwo\;|\;t\z1s+s\ztwos=st\,\},$$ and $$N=
\{\,(e^{2\pi iX},e^{2\pi i\frac{s}{t}X})\;|\; X\in\R\,\}.$$ The space $\squo$ is
covered by two charts $$\us=\left\{\,[z_1:z_2]\in \squo\;|\;z_2\neq
0\,\right\}\quad\hbox{and}\quad \un=\left\{\,[z_1:z_2] \in \squo\;|\;z_1\neq
0\,\right\}.$$The corresponding local models are defined by
$\uts=\{z\in\C\;|\;|z|<\sqrt{s}\}\subset\vts$ acted on by $\G_{\sc s}$ and by
$\utn=\{z\in\C\;|\;|z|<\sqrt{t}\}\subset\vtn$ acted on by $\G_{\sc n}$. To complete
the  picture we need the homeomorphisms: $$\begin{array}{cccc} \Phi_{\sc
s}\,\colon\,&\uts/\G_{\sc s}&\longrightarrow &\us\\
&[z]&\longmapsto & \left[z:\sqrt{t-\frac{t}{s}|z|^2}\right]\end{array},\quad
\begin{array}{cccc}
\Phi_{\sc n}\,\colon\,&\utn/\G_{\sc n}&\longrightarrow &\un\\
&[w]&\longmapsto & \left[\sqrt{s-\frac{s}{t}|w|^2}:w\right]\end{array}.$$ Now we are
ready to write the local lifts $\chis$ and $\chin$ of the diffeomorphism $\chi$.
Notice that $$(\chi\circ\Phi_{\sc s})([z])=
\left[z:\sqrt{t-\frac{t}{s}|z|^2}\right]_{\Nc}=
\left[z(t-\frac{t}{s}|z|^2)^{-\frac{s}{2t}}:1\right]_{\Nc}.$$ Therefore a local lift
$\chis:\uts\lorw\vts$ is given by the equivariant diffeomorphism $$
\chis(z)=z\left(t-\frac{t}{s}|z|^2\right)^{-\frac{s}{2t}} $$ Analogously a local lift
$\chin:\utn\lorw\vtn$ is given by the equivariant diffeomorphism $$
\chin(w)=w\left(s-\frac{s}{t}|w|^2\right)^{-\frac{t}{2s}}. $$ We exhibit now the
complex quotient as a K\"ahler quasifold. We define the K\"ahler form $\omega$ by
giving its local lifts on $\vts$ and $\vtn$:
\begin{equation}\label{kahlersfera}
\tilde{\omega}_{\sc s}=\frac{1}{2\pi i}\frac{s}{t}\frac{1}
{\left(\frac{1}{s}+\frac{1}{t}|z|^{(1+\frac{s}{t})}\right)^2}\, dz\wedge
d\overline{z},\quad\quad \tilde{\omega}_{\sc n}=\frac{1}{2\pi i}\frac{t}{s}\frac{1}
{\left(\frac{1}{t}+\frac{1}{s}|w|^{(1+\frac{t}{s})}\right)^2}\, dw\wedge
d\overline{w}.
\end{equation}
They are invariant under the action of $\G_{\sc s}$ and $\G_{\sc n}$ respectively and
it is a straightforward computation to check that $\chis^*\tilde{\omega}_{\sc
s}=\frac{1}{2\pi i}\,dz\wedge d\overline{z}$ and $\chin^*\tilde{\omega}_{\sc
n}=\frac{1}{2\pi i}\,dw\wedge d\overline{w}$. On the other hand these are precisely
the local lifts of the symplectic form on $\uts$ and $\utn$ respectively, so $\omega$
is the pullback via $\chi^{-1}$ of the symplectic form on $\squo$ and hence defines a
K\"ahler form on the complex quotient.

To conclude the discussion of this significant example, it is worthwhile to check
directly that the local forms given in (\ref{kahlersfera}) fulfill the definition of
differential form. We have to show that they behave correctly under the change of
charts. The first move is to pull them back to $\wts$ and to $\wtn$, which are both
equal to $\C$. We obtain $$\omega^{\#}_{\sc s}= \frac{1}{2\pi
i}\frac{s}{t}\frac{e^{(\zeta+\overline{\zeta})} }
{\left(\frac{1}{s}+\frac{1}{t}e^{(\zeta+\overline{\zeta})\frac{1}{2} (1+\frac{s}{t})}
\right)^2}\,d\zeta\wedge d\overline{\zeta},\quad\quad \omega^{\#}_{\sc
n}=\frac{1}{2\pi i}\frac{t}{s} \frac{1} {\left(\frac{1}{t}+\frac{1}{s}
e^{(\eta+\overline{\eta})\frac{1}{2}(1+\frac{t}{s})} \right)^2}\,d\eta\wedge
d\overline{\eta}$$ which are invariant under the respective action of $\Lambda_{\sc
s}$, $\Lambda_{\sc n}$. It is easy to check that the mapping
$\zeta\mapsto\eta=-\frac{s}{t}\zeta$ pulls $\omega^{\#}_{\sc n}$ back to
$\omega^{\#}_{\sc s}$.}
\end{ex}

\noindent \small{\sc Dipartimento di Matematica Applicata "G. Sansone", Via S. Marta
3, 50139 Firenze, ITALY,
                         {\tt mailto:fiamma@dma.unifi.it}\\
                        and\\
                        Laboratoire Dieudonn\'e,
                        Universit\'e de Nice, Parc Valrose, 06108 Nice
                        Cedex 2, FRANCE, {\tt mailto:elisa@alum.mit.edu}}

\begin{thebibliography}{AA}

\bibitem[D]{d} T. Delzant, Hamiltoniens p\'eriodiques et image convexe
de l'application moment, {\em Bull. S.M.F.} {\bf 116} (1988), 315--339.

\bibitem[G]{vg} V. Guillemin, Moment maps and combinatorial invariants
of Hamiltonian $T^n$-spaces, Progress in Mathematics 122, Birkh\"auser, Boston, 1994.

\bibitem[KN]{kempfness} G. Kempf and L. Ness, The length of vectors in
representation spaces, Algebraic Geometry, Proceedings, Copenhagen, 1978, Lecture
Notes in Math. {\bf 732} (1979), 233--244.

\bibitem[LT]{lt} E. Lerman and S. Tolman, Hamiltonian torus actions
on symplectic orbifolds and toric varieties, {\em Trans. A.M.S.} {\bf 349} No. 10
(1997), 4201--4230.

\bibitem[K]{kirwan} F. Kirwan, Cohomology of quotients in symplectic and algebraic
geometry, Mathematical Notes 31, Princeton University Press, 1984.

\bibitem[N]{ness} L. Ness, A stratification of the null cone via the moment
map, {\it Amer. J. Math.} {\bf 106} (1984), 1281--1329.

\bibitem[P]{p} E. Prato, Simple non-rational convex polytopes via
symplectic geometry, preprint {\tt math.SG/9904179}, to appear in {\em Topology}.

\end{thebibliography}
\end{document}